\documentclass[12pt, reqno]{amsart}

\usepackage{amsfonts,amsmath,amsthm,amssymb}
\usepackage{latexsym}

\begin{document}

{\theoremstyle{plain}
  \newtheorem{theorem}{Theorem}[section]
  \newtheorem{corollary}[theorem]{Corollary}
  \newtheorem{proposition}[theorem]{Proposition}
  \newtheorem{lemma}[theorem]{Lemma}
  \newtheorem{question}[theorem]{Question}
  \newtheorem{conjecture}[theorem]{Conjecture}
  \newtheorem{claim}[theorem]{Claim}
}

{\theoremstyle{definition}
  \newtheorem{definition}[theorem]{Definition}
  \newtheorem{remark}[theorem]{Remark}
  \newtheorem{example}[theorem]{Example}

}

\title [Formal prime Ideals of infinite value]{Formal prime ideals of infinite value and their algebraic resolution}

\author{Steven Dale Cutkosky}
\address{Steven Dale Cutkosky, Department of Mathematics, University of Missouri, Columbia, MO 65211, USA}
\email{cutkoskys@missouri.edu}
\thanks{The first author was partially supported by NSF}

\author{Samar ElHitti}
\address{Samar ElHitti, Department of Mathematics, New York City College of Technology, 300 Jay street, Brooklyn, NY 11201, USA}
\email{selhitti@citytech.cuny.edu}
\urladdr{http://websupport1.citytech.cuny.edu/faculty/selhitti}

\thanks{The second author was partially supported by the PSC-CUNY Award \# 60070-39 40}

\begin{abstract}
Suppose that $R$ is a local domain essentially of finite type over a field of characteristic $0$, and $\nu$ a valuation of the quotient field of $R$ which dominates $R$. The rank of such a valuation often increases upon extending the valuation to a valuation dominating $\hat R$, the completion of $R$. When the rank of $\nu$ is $1$, Cutkosky and Ghezzi handle this phenomenon by resolving the prime ideal of infinite value, but give an example showing that when the rank is greater than $1$, there is no natural ideal in $\hat R$ that leads to this obstruction. We extend their result on the resolution of prime ideals of infinite value to valuations of arbitrary rank.
\end{abstract}



\maketitle

\section{Introduction}

Since Zariski introduced general valuation theory into algebraic geometry,
valuations have been important in addressing  problems on resolution of singularities.

Suppose that $K$ is an algebraic function field over a base field $k$, and $V$ is a
valuation ring of $K$. $V$ determines a unique center on a proper variety $X$ whose
function field is $K$. The valuation gives a way of reducing a global problem on
$X$, such as resolution, to a local problem, studying the local rings of centers of
$V$ on different varieties $X$ whose function field is $K$.

The valuation theoretic analogue of resolution of singularities is local
uniformization.
The problem of local uniformization is to find, for a fixed valuation ring $V$ of a function field $K$ over $k$,  a regular local ring $R$ essentially
of finite type over $k$ with quotient field $K$ such that  $V$
dominates $R$. That is, $R\subset V$ and $m_V\cap R=m_R$. In 1944, Zariski
\cite{OZ2} proved local uniformization over fields of characteristic zero.
A consequence is the following theorem.

\begin{theorem}(Zariski) Suppose that $R$ is a regular local ring
which is essentially of finite type over a field of characteristic zero, which
is dominated by a valuation ring $V$.
Suppose that $f\in R$. Then there exists a birational extension of regular local
rings $R \rightarrow R_1$ such that $R_1$ is dominated by $V$, and
ord$_{R_1}\overline f \leq 1$ where $\overline f$ is the strict transform of $f$ in
$R_1$. If $V$ has rank $1$, then there exists $R_1$ such that $\overline f$ is a
unit in $R_1$.
\end{theorem}

Zariski first proved local uniformization for two-dimensional function fields over
an algebraically closed field of characteristic zero in \cite{OZ1}. He later proved
local uniformization for algebraic function fields of characteristic $0$ in
\cite{OZ2}. Hironaka proved resolution of singularities of characteristic zero
varieties in 1964 \cite{Hir}.

Abhyankar has proven local uniformization and resolution of singularities in positive characteristic for
two dimensional function fields, surfaces and three dimensional varieties (characteristic $p>5$)
\cite{Ab1}, \cite{Ab3}.

Recently, there has been progress on local uniformization in positive characteristic, including the work of Cossart and Piltant \cite{CP1}, \cite{CP2}, Kuhlmann \cite{Ku1}, Knaf and Kuhlmann \cite{KK}, Spivakovsky \cite{Sp}, \cite{S}, Temkin \cite{Te} and  Teissier \cite{T}. Some recent progress on understanding valuations
in the context of algebraic geometry has been made by
Favre and Jonnson \cite{FJ}, Ghezzi, H\`a and Kascheyeva \cite{GHK}, Vaqui\'e \cite{Va} and others.

One of the most important techniques in studying resolution problems is to pass to
the  completion of a local ring (the germ of a singularity). This allows us to
reduce local questions on singularities to problems on power series.

Let $R$ be a regular local ring, and let $V$ be a valuation ring of the quotient field of $R$ which dominates $R$. Let $\nu$ be a valuation whose valuation ring
is $V$.

A question which arises on completion is if the following generalization of
local uniformization is true:

\begin{question}\label{question1}
Given a reduced element $f \in \hat R$, when does there exist a birational extension $R\rightarrow R_1$
where $R_1$ is a regular local ring dominated by $V$ such that $ord_{\hat
R_1}\overline f \leq 1$, where $\overline f$ is a strict transform of $f$ in $\hat
R_1$?
\end{question}

Question \ref{question1} has a positive answer if $R$ is a regular local ring of dimension 2, since a germ of a curve singularity can be resolved by blowing up points.  If $R$ is essentially of finite type over a field of characteristic zero, of arbitrary dimension, and $f\notin Q(\hat R)$, then a positive answer to Question \ref{question1} is a consequence of Theorem \ref{res}.

The answer to question \ref{question1} is however generally no. We give a simple  example stated below,  which comes from a discrete valuation.

\begin{example}\label{Example1}
There exists a discrete valuation ring $V$ dominating a regular local ring $R$ of
dimension $3$ such that for all integers $r\geq 2$, there exists an irreducible element $f \in \hat R$ such that for
all birational extensions $R\rightarrow R_1$ of regular local rings dominated by $V$, a strict transform of $f$ has order $\geq r$ in $\hat R_1$.
\end{example}

The example is constructed in Section \ref{example}.
It can be understood in terms of an extension of our valuation ring $V$ to
a valuation ring dominating $\hat R$. It is a fact that the rank of a valuation $V$ dominating $R$ often increases when extending the valuation to a valuation ring $\hat V$ dominating $\hat R$. In the comment after the statement of Theorem C on page 177 of \cite{HS}, it is shown that
if the extension of $V$ to a valuation ring of the quotient field of $\hat R$
which dominates $\hat R$ is not unique, then $V$ has extensions of at least two ranks. Some papers where
this phenomonen is studied are Spivakovsky \cite{S}, Heinzer and Sally \cite{HS}, and Cutkosky and Ghezzi \cite{CG}.

The first measure of this increase of rank is the prime ideal
$Q(\hat R)$, defined in Section \ref{introduction}. This ideal is known as the prime ideal of infinite value. $Q(\hat R)$ is generated by Cauchy sequences $\{f_n\}$ in $R$ whose values are eventually larger than
$t\nu(m_R)$ for any multiple $t$  of $\nu(m_R)$.

This prime has been previously defined and studied by Teissier \cite{T}, Cutkosky \cite{C}, Cutkosky and Ghezzi \cite{CG} and Spivakovsky \cite{S}. We show in Theorem \ref{res}, stated below, that if $f\not\in Q(\hat R)$, then Question \ref{question1} does have a positive answer.

\begin{theorem}\label{res} Suppose that $R$ is a regular local ring which is essentially of finite type over a field of characteristic zero, and is dominated by a valuation ring $V$. Suppose that $f\in\hat R$ is such that   $f\not\in
Q(\hat R)$. Then there exists a sequence of monoidal transforms
$R\rightarrow R_1$ along $V$ such that a strict transform $\overline f$ of $f$ in $\hat R_1$ is a unit.
\end{theorem}

Theorem \ref{res} is proven in Section \ref{example}.

If the rank of $V$ is 1, then there is a unique extension of $V$ to the quotient field of $\hat R/Q(\hat R)$ dominating $\hat R/Q(\hat R)$, and the rank of this extension does not increase. The prime ideal $Q(\hat R)$ led to this obstruction.

In spite of the fact that we cannot always
resolve the singularity of $f=0$ by a birational extension of $R$, we can resolve
the formal singularity, whose local ring is $\hat R/Q(\hat R)$, by a birational
extension of $R$. This is proven for valuations of rank $1$ by Cutkosky and Ghezzi
\cite{CG}.

\begin{theorem}\label{CG}(Cutkosky, Ghezzi)
Suppose that $R$ is a local domain which is essentially of finite type over a field of characteristic zero and $V$ is a rank $1$ valuation ring which dominates $R$.  Then there exists a birational extension $R\rightarrow R_1$
where $R_1$ is a regular local ring dominated by $V$ such that $Q(\hat R_1)$ is a
regular prime.
\end{theorem}

A regular prime is a prime ideal $P$ in a regular local ring $R$ such that $R/P$ is also a regular local ring.

In Example 4.3 \cite{CG} Cutkosky and Ghezzi consider a valuation $\nu$ of rank $2$ and give two different extensions of $V$ to $\hat R$ one of rank $2$ and the other of rank $3$.
Thus when rank of $V$ is greater than $1$, there is no natural ideal in $\hat R$ that obstructs the jumping of the rank of an extension of $V$ to $\hat R$. This obstruction will be obtained in a series of prime ideals in quotient rings in $\hat R$ as stated in section \ref{mainresult}.

The main result of this paper is to generalize Cutkosky and Ghezzi's result on the resolution of prime ideals of infinite value to valuations of arbitrary rank. We state our result below.

\begin{theorem}\label{SEH}

 Suppose that $R$ is a local domain which is essentially
of finite type over a field $k$ of characteristic zero, and $V$ is a
 valuation ring of the quotient field $K$ of $R$ which dominates $R$. Let

 $$
 (0)=P_{V,t}\subset\cdots\subset P_{V,1}\subset P_{V,0}=m_V
 $$
\vskip .2truein
\noindent  be the chain of prime ideals of $V$.
Then there exists a sequence of monoidal transforms
 $R\rightarrow R_1$ such that $R_1$ is a
regular local ring
and $V$ dominates $R_1$. Further, $P_{R_1,i}=P_{V,i}\cap R_1$ are distinct regular primes such that $(V/P_{V,i})_{P_{V,i}}$ is algebraic over $(R/P_{R_1,i})_{P_{R_1,i}}$ for all
$i$, and

$$
Q(\widehat{(R_1)_{P_{R_1,i}}})\subset \widehat{(R_1)_{P_{R_1,i}}}
$$
\noindent  are regular primes for all $i$.
\end{theorem}

Theorem \ref{SEH} is stated and proved in section \ref{mainresult}.
A different proof of Theorem \ref{SEH} using a generalization of Zariski's Perron transforms to higher rank valuations is given in ElHitti's Ph.D. thesis \cite{ElH}.

The following related question is raised by Teissier and Spivakovsky in their work on local uniformization.

\begin{question}\label{questionlast} Suppose that $R$ is a local domain, containing a field $k$
(of any characteristic) which is dominated by $V$. Does there exist a birational extension $R\rightarrow R_1$ such that $R_1$ is regular
and a  prime ideal $P\subset\widehat{R_1}$  such that $P\cap R_1=(0)$, $\widehat{R_1}/P$ is regular and
there is a unique extension of $V$ to the quotient field of $\widehat{R_1}/P$
which dominates $\widehat{R_1}/P$ and the rank does not increase?
\end{question}

In the case when $V$ has rank $1$, and characteristic zero, a positive answer to
Question \ref{questionlast} follows from Cutkosky and Ghezzi \cite{CG}.

\section{Notations and Preliminaries}\label{introduction}
The maximal ideal of a quasi local ring $R$ will be written $m_R$. Suppose that $R\subset S$ is an inclusion of quasi local rings. We will say that $R$ dominates $S$ if $m_S\cap R =m_R$.
If $R$ is a local ring (a noetherian quasi local ring), $\hat R$ will denote the completion of $R$ at its maximal ideal. A prime ideal $P$ is a regular prime if $R/P$ is a regular local ring.

Good introductions to valuation theory can be found in Chapter VI of \cite{ZS} and in \cite{Ab4}.
We will summarize a few basic results and set up the notation which we will use.

Suppose that $V$ is a valuation of rank $t>1$, with quotient field $K$ and valuation ring $\Gamma=\Gamma_V$.
Let $\nu$ be a valuation of $K$ whose valuation ring is $V$. Suppose that $V$ dominates a noetherian local
domain $R$ whose quotient field is $K$. Then the rank $t$ of $V$ is finite (\cite{Ab2} or Appendix 2 of \cite{ZS}).

Let
$$
(0)=P_{V,t}\subset \cdots \subset P_{V,1}\subset P_{V,0}=m_V
$$
be the chain of distinct prime ideals in $V$.
There is a one-to-one order reversing correspondence between the isolated subgroups $\Gamma_i$ of $\Gamma_V$ and the the prime ideals of $V$ (cf. Sections 8, 9, 10 of \cite{Ab4} and chapter VI, section 10 of \cite{ZS}), giving the sequence
$$
\{0\}=\Gamma_0\subset \Gamma_1\subset \cdots\subset \Gamma_t=\Gamma_V
$$
of isolated subgroups of $\Gamma_V$. For $0\le i\le t$, let
$$
U_i=\{\nu(a)\mid a\in P_{V,i}\}.
$$
$\Gamma_i$ is defined to be the complement of $U_i$ and $-U_i$ in $\Gamma_V$.

For $i>j$, $\nu$ induces a valuation on the field $(V/P_{V,j})_{P_{V,j}}$ with valuation ring $(V/P_{V,j})_{P_{V,i}}$ and value group $\Gamma_j/\Gamma_i$. In particular if $j=i+1$, $\Gamma_j/\Gamma_i$ has rank $1$.

Let
\begin{equation}\label{eq1}
(0)=P_{R,t}\subset \cdots \subset P_{R,1}\subset P_{R,0}=m_{R}
\end{equation}
 be the induced chain of prime ideals in $R$, where $P_{R,i}=P_{V,i}\cap R$.\label{inducedchain}
For all $i$, $V_{P_{V,i}}$ is a valuation ring of $K$ dominating $R_{P_{R,i}}$. Let $\nu_i$ be the valuation associated to  $V_{P_{V,i}}$.

Suppose that $f\in \widehat{R_{P_{R,i}}}$.
Consider the following condition
 on a Cauchy sequence $\{f_n\}$ of elements $f_n\in R_{P_{R,i}}$ converging to $f$:
\begin{equation}\label{eq2}
\text{For all $l\in \mathbb{N}$, there exists $n_l\in \mathbb{N}$ such that
$\nu_i(f_n)\ge l\nu(m_{R_{P_{R,i}}})$ if $n\ge n_l$.}
\end{equation}

The condition (\ref{eq2}) is independent of Cauchy sequence $\{f_n\}$ converging to $f$, although the
specific numbers $n_l$ depend on the Cauchy sequence.

Define the prime ideal $Q(\widehat{R_{P_{R,i}}})\subset \widehat{R_{P_{R,i}}}$ for the
valuation ring $V_{P_{V,i}}$
 by
$$
Q(\widehat{R_{P_{R,i}}})=\left\{
\begin{array}{ll}
f\in {\widehat{R_{P_{R,i}}}}\mid
\text{ $f$ satisfies (\ref{eq2})}
\end{array}\right\}.
$$
We have
$$
Q(\widehat{R_{P_{R,i}}})\cap R_{P_{R,i}}=P_{R,i+1}R_{P_{R,i}}
$$
for all $i$.

If $i=0$, then we have $\nu_0=\nu$, $R_{P_{R,0}}=R$, $m_{R_{P_{R,0}}}=m_R$, and
$$
Q(\hat R)=\{f\in\hat R\mid f\mbox{ satisfies (\ref{eq2})}\}.
$$
We have that  $Q(\hat R)\cap R= P_{R,1}$, which is the zero ideal if $V$ has rank $1$.

Suppose that $R$ is a regular local ring. A monoidal transform $R\rightarrow R_1$ is a birational extension of local domains such that $R_1=R[\frac{P}{x}]_m$ where $P$ is a regular prime ideal of $R$, $0\neq x \in P$ and $m$ is a prime ideal of $R\left[\frac{P}{x}\right]$ such that $m\cap R=m_R$. $R\rightarrow R_1$ is called a quadratic transform if $P=m_R$.

If  $R\rightarrow R_1$ is a monoidal transform, then there exists a regular system of parameters $(x_1, \dots, x_n)$ in $R$ and $r\leq n$ such that
$$
R_1=R\left[\frac{x_2}{x_1},\dots, \frac{x_r}{x_1}\right]_m.
$$
If $V$ is a valuation ring dominating $R$, then $R\rightarrow R_1$ is called a monoidal transform
along $V$ if $V$ dominates $R_1$.

Suppose that $I\subset R$ is an ideal. We define the strict transform $I_1$ of $I$ in $R_1$ by
$$
I_1=\bigcup_{j=1}^{\infty} (IR_1:P^jR_1).
$$

For an ideal $J\subset \hat R$, we define the strict transform $J_1$ of $J$ in $\hat R_1$ by
$$
J_1=\bigcup_{j=1}^{\infty}(J\hat R_1:P^j\hat R_1).
$$

If $f\in R$, a strict transform $\overline f$ of $f$ in $R_1$ is a generator of the
strict transform in $R_1$ of the ideal generated by $f$ in $R$. A strict transform of an element in $\hat R$ is  defined in the same way.

We will make use of results on resolution of singularities by Hironaka \cite{Hir}, in a form suitable to
application to monoidal transforms along valuations, from Chapter 2 of  \cite{C}.

Suppose that $R$ is a local domain, $P\subset R$ is a prime ideal, and $f\in R$. We define $\mbox{ord}_P(f)$
to be the largest integer $n$ such that $f\in P^{(n)}$. Here $P^{(n)}$ denotes the $n$-th symbolic power of $P$.
If $P$ is  a regular prime in a regular local ring $R$, then $P^{(n)}=P^n$ for all $n$; that is, the ordinary and symbolic powers agree.
We will write $\mbox{ord}_R(f)=\mbox{ord}_{m_R}(f)$. If $\gamma\subset \mbox{Spec}(R)$ is an integral subvariety, and $I_{\gamma}$ is the prime ideal of $\gamma$, then we define $\mbox{ord}_{\gamma}(f)=\mbox{ord}_{I_{\gamma}}(f)$.

\section{Algebraic resolution of  formal series}\label{example}

In this section we give proofs of Theorem \ref{res} and Example \ref{Example1}
stated in the introduction. We use the notation established in Section \ref{introduction}.
\vskip .2truein

\begin{proof}[\bf Proof of Theorem \ref{res}]
Let $\nu$ be a valuation whose valuation ring is $V$. By Theorem 2.9 \cite{C}, there exists a sequence of monoidal transforms $R\rightarrow R_1$ along $V$ such that the strict transform of $P_{R,1}$ in $R_1$ is a regular prime, which is thus necessarily $P_{R_1,1}$. Set $A_1=R_1/P_{R_1,1}$. $\nu$ induces a rank $1$ valuation $\overline \nu$, on the quotient field $K$ of $A_1$, which dominates $A_1$ and has valuation ring $\overline V=(V/P_{V,1})\cap K$. Let $\tilde f$ be the residue of $f$ in $\widehat{A_1}$. Let $s$ be the rational rank of $\overline\nu$.
We have that $f\notin Q(\widehat{A_1})$ (computed with respect to $\overline\nu$), since
$f\not\in Q(\widehat{R_1})$. Thus  $\overline\nu(\tilde f)<\infty$ (in the notation of \cite{C}). By Theorem 4.8 \cite{C}, there exists a sequence of monoidal transforms
$$
A_1\rightarrow \cdots \rightarrow A_a
$$
along $\overline V$ such that
$$
\tilde f=x_1^{a_1}\cdots x_s^{a_s}u+h
$$
where $x_1,\ldots,x_s,x_{s+1},\ldots,x_n$ are a  regular system of parameters in $A_a$, $x_1\cdots x_s=0$ is a local equation of the exceptional divisor of $\mbox{Spec}(A_a)\rightarrow \mbox{Spec}(A_1)$, $u\in\widehat{A_a}$ is a unit series, $a_1,\ldots,a_s$ are positive integers and $h\in m_{\widehat{A_{a}}}^N$ where $N\overline\nu(m_{A_{a}})
>\overline\nu(x_1^{a_1}\cdots x_s^{a_s})$.

Now by (2) of Theorem 4.10 \cite{C}, there exists a sequence of monoidal transforms along $\overline V$
$$
A_a\rightarrow \cdots \rightarrow A_{a+b}
$$
and regular parameters $y_1,\ldots, y_s,y_{s+1},\ldots, y_n$ in $A_{a+b}$
such that $y_1\cdots y_s=0$ is a local equation of the exceptional divisor of $\mbox{Spec}(A_{a+b})\rightarrow \mbox{Spec}(A_1)$ and
$$
\tilde f=y_1^{b_1}\cdots y_s^{b_s}\tilde u
$$
 where $b_1,\ldots,b_s$ are positive integers and $\tilde u\in \widehat{A_{a+b}}$ is a unit series.

We may now construct a sequence of monoidal transforms
$$
R_1\rightarrow \cdots \rightarrow R_{a+b}
$$
along $V$ (as in the last part of the proof of Theorem \ref{SEH} in the next section) such that
$$
(R_j)_{P_{R_j,1}}\cong (R_1)_{P_{R_1,1}}\mbox{ and }R_j/P_{R_j,1}\cong A_j
$$
for $1\le j\le a+b$.

 There exists a regular system of parameters $z_1,\ldots,z_c$ in $R_{a+b}$ such that the residue of $z_i$ in $A_{a+b}$ is $y_i$ for $1\le i\le n$ and $P_{R_{a+b},1}=(z_{n+1},\ldots,z_c)$.
Let
$$
B=R_{a+b}\left[\frac{z_{n+1}}{z_1^{b_1}\cdots z_s^{b_s}},\ldots,
\frac{z_c}{z_1^{b_1}\cdots z_s^{b_s}}\right].
$$
We have $B\subset V$. let $R'=B_{B\cap m_V}$. $R_{a+b}\rightarrow R'$
factors as a sequence of monoidal transforms along $V$. Define a regular system of parameters $z_1(1),\ldots, z_c(1)$ in $R'$ by
$$
z_i=\left\{\begin{array}{ll}
z_i(1)&\mbox{ if }1\le i\le n\\
z_1(1)^{b_1}\cdots z_s(1)^{b_s}z_i(1)&\mbox{ if }n+1\le i\le c.
\end{array}\right.
$$
There exists a unit series $\Lambda\in\widehat{R_{a+b}}$ and $g\in P_{R_{a+b},1}\widehat{R_{a+b}}$ such that
$$
f=z_1^{b_1}\cdots z_s^{b_s}\Lambda+g.
$$
Thus
$$
f=z_1(1)^{b_1}\cdots z_s(1)^{b_s}(\Lambda+g')
$$
where $g'\in P_{R',1}\widehat{R'}$. The unit series $\Lambda+g'$ is a strict transform $\overline f$ of $f$ in $\widehat{R'}$.
\end{proof}

Theorem \ref{res} does not generalize if $f\in Q(\hat R)$, as is shown by  Example \ref{Example1}, stated in the introduction.
\vskip .2truein

\begin{proof}[\bf Proof of Example \ref{Example1}]
Let $k$ be a field, $k[[t]]$ be a power series ring over $k$ and $t, \phi(t), \psi(t) \in k[[t]]$ be algebraically independent elements of positive order.
We have an inclusion $k(x,y,z)\hookrightarrow k((t))$ of $k$-algebras defined by the substitutions $x=t, y=\phi(t)$ and $z=\psi(t)$.
The order valuation on $k((t))$ (with valuation ring $k[[t]]$) restricts to a discrete rank $1$ valuation $\nu$ on $k(x,y,z)$, dominating $R=k[x,y,z]_{(x,y,z)}$.

The kernel $Q(\hat R)=(y-\phi(x),z-\psi(x)) \subset \hat R$  of the induced homomorphism $\hat R\rightarrow k[[t]]$ is a regular prime of height $2$, and it defines a nonsingular curve $\gamma \subset \text{Spec }(\hat R)$, with ideal $I_{\gamma}=Q(\hat R)$.
We have $Q(\hat R)\cap R=(0)$.

Suppose that $r \in \mathbb{N}$ (with $r \geq 2$).
Let
$$
f=(y-\phi(x))^{r}+(z-\psi(x))^{r+1} \in \hat R.
$$
We have
$$
\mbox{ord}_{\hat R}(f)=r\mbox{ and }\mbox{ord}_{\gamma}(f)=\mbox{ord}_{(y-\phi(x),z-\psi(x))}(f)=r.
$$
Suppose that $R\rightarrow R_1$ is a birational extension where $R_1$ is a regular local ring dominated by $V$.

Let $\overline f$ be the strict transform of $f$ in $R_1$. The ideal $Q(\widehat{R_1})$ is the kernel of
the induced homomorphism $\widehat{R_1}\rightarrow k[[t]]$. Let $S_1=R_1\otimes_R\hat R$, and let
$\gamma'$ be the strict transform of $\gamma$ in $\mbox{Spec}(S_1)$, with ideal sheaf $I_{\gamma'}$.
 Since $I_{\gamma}\cap R=\{0\}$ and $R_1$ is a local ring of the blow up of a nonzero  ideal in $R$, $(S_1)_{I_{\gamma'}}\cong (\hat R)_{I_{\gamma}}$. Thus
$$
r=\mbox{ord}_{\gamma}(f)=\mbox{ord}_{\gamma'}(f)=\mbox{ord}_{\gamma'}(\overline f).
$$

From the natural inclusions $\hat R/I_{\gamma}\subset S_1/I_{\gamma'}\subset k[[t]]$, and the fact that $\hat R/I_{\gamma}\cong k[[t]]$, we have $S_1/I_{\gamma'}\cong k[[t]]$. Thus $I_{\gamma'}$ is a regular prime in the regular local ring $S_1$, and thus since $\widehat{S_1}\cong \widehat{R_1}$ and $I_{\gamma'}\widehat{R_1}=Q(\widehat{R_1})$, we have
$$
\mbox{ord}_{\widehat{R_1}}(\overline f)\ge \mbox{ord}_{Q(\widehat{R_1})}(\overline f)=\mbox{ord}_{\gamma'}(\overline f)=r.
$$
\end{proof}

In Example \ref{Example1}, the rank of the valuation  $\nu$ must increase when passing to the completion $\hat R$, since $Q(\hat R)\ne 0$.
The following is a construction of a valuation $\hat {\nu}$ of the quotient field of $\hat{R}$ which extends $\nu$ and dominates $\hat R$.
For a nonzero element $\delta\in \hat R$, write
$
\delta=(y-\phi(x))^{\alpha}g(x,y,z)
$
where $g(x,\phi(x),z)\neq 0$.
Write $g(x,\phi(x),z)=(z-\psi(x))^{\beta}h(x,\phi(x),z)$ where $h(x,\phi(x),\psi(x))\neq 0$. Set $\gamma = \mbox{ord}( h(t,\phi(t),\psi(t)))$.

Define
$$\hat {\nu}(\delta)=(\alpha,\beta,\gamma)\in \mathbb{Z}^3$$
where $\mathbb{Z}^3$ has the lexicographic order. $\hat {\nu}$ extends to a rank $3$ valuation of the quotient field of $\hat{R}$ which  dominates $\hat R$ and extends $\nu$.

\section{A resolution theorem for formal ideals along a high rank valuation}\label{mainresult}

In this section we prove our main resolution theorem, Theorem \ref{SEH},
which is stated in the introduction. We use the notation established in Section \ref{introduction}.

\begin{proof}[\bf Proof of Theorem \ref{SEH}]
Let $\nu$ be a valuation of the quotient field of $R$ whose valuation ring is $V$.
$R$ is a local domain which is essentially of finite type over $k$. Thus there exists a
 regular local ring $T$ which is essentially of finite type over $k$, and a prime ideal $P$ in $T$ such that
$R=T/P$.
Let $\nu_1$ be the $PT_P$-adic valuation of the regular local ring $T_P$, and let $w=\nu \circ \nu_1$
be the composite valuation (Section 10 of \cite{Ab4} or Section 11, Chapter VI of \cite{ZS}). $w$ is a rank $t+1$ valuation that dominates $T$. Let $W$ be the valuation ring of $w$, and
$$
(0)=P_{W,t+1}\subset \cdots \subset P_{W,0}=m_W
$$
be the chain of prime ideals in $W$. We have that $V/P_{V,i}\cong W/P_{W,i}$ for $0\leq i\leq t$.

By Theorem 2.9 \cite{C} there exists a sequence of monoidal transforms $T\rightarrow T_0$ along $w$ such that $T_0$ is a regular local ring, the strict transform $P_0$ of $P$ in $T_0$ is a regular prime, and
 $R\rightarrow R_0=T_0/P_0$ factors as a sequence of monoidal transforms along $V$.
We will first show that there exists a a sequence of monoidal transforms
$$
R \rightarrow R_1 \rightarrow \dots \rightarrow R_t
$$
along $V$ such that
\begin{equation}\label{algebraic}
(V/P_{V,i})_{P_{V,i}}\mbox{ is algebraic over } (R_t/P_{R_t,i})_{P_{R_t,i}}\mbox{ for all $i$}.
\end{equation}
and
\begin{equation}\label{regdistinct}
P_{R_t,i}\mbox{ are regular and distinct prime ideals in }R_t\mbox{ for all $i$}.
\end{equation}

Let
$$
(0)=P_{R_0,t}\subset \cdots \subset P_{R_0,1}\subset P_{R_0,0}=m_{R_0}
$$
 be the induced chain of prime ideals in $R_0$ as defined in (\ref{eq1}).
We have that $V_{P_{V,i}}$ dominates $(R_0)_{P_{R_0,i}}$.

Let $d_i= \text{trdeg} _{(R_0)_{P_{R_0,i}}/P_{R_0,i}(R_0)_{P_{R_0,i}}}V_{P_{V,i}}/P_{V,i} V_{P_{V,i}}$.
By Abhyankar's inequality, $d_i$ is finite (Theorem 1 \cite{Ab2}). For all $i$, let $\bar{t}_{i1},\dots ,\bar{t}_{id_i}$ be a transcendence basis of $(V/P_{V,i})_{P_{V,i}}$
over $(R_0/P_{R_0,i})_{P_{R_0,i}}$.
For all $i$, let $t_{i1}, \dots, t_{id_i}$ be a lift of this transcendence basis to $V$, after possibly replacing some $t_{ij}$ with $1/t_{ij}$, so that we have  $\nu(t_{ij})\geq 0$ for $j=1, \dots, d_i$ and $i=0, \dots, t$.

By Theorem 2.7 \cite{C} there exists a sequence of monoidal transforms $R_0\rightarrow R_1$ along $V$ such that $t_{i1}, \dots, t_{id_i} \in R_1$ for all $i$. Thus

\begin{equation}\label{thiscondition}
(V/P_{V,i})_{P_{V,i}} \text{ is algebraic over } (R_1/P_{R_1,i})_{P_{R_1,i}}  \text{for all } i.
\end{equation}

Consider the chain of  prime ideals in $R_1$ given by $P_{R_1,i}=P_{V,i}\cap R_1$ for  $i=0, \dots , t$.
Suppose that $P_{R_1,i} = P_{R_1,i-1}$ for some $i$.
Choose $f\in P_{V,i-1}-P_{V,i}$. Let $S=R_1[f]$. Define $P_{S,j}=P_{V,j}\cap S$ for $0\le j\le t$.
we have an inequality of heights
$$
\mbox{ht}(P_{S,i})<\mbox{ht}(P_{S,i-1}),\mbox{ since } f\in P_{S,i-1}-P_{S,i}.
$$
For all $j$, we have that $(S/P_{S,j})_{P_{S,j}}$ is an intermediate field between $(R_1/P_{R_1,j})_{P_{R_1,j}}$ and $(V/P_{V,j})_{P_{V,j}}$. Thus $(S/P_{S,j})_{P_{S,j}}$ is algebraic over $(R_1/P_{R_1,j})_{P_{R_1,j}}$ for all $j$.
Since $S$ is a birational extension of $R_1$, we now conclude  that $\mbox{ht}(P_{R_1,j})=\mbox{ht}(P_{S,j})$ for all $j$ by the dimension formula (Theorem 15.6 \cite{M}). This contradicts our assumption that $P_{R_1,i}=P_{R_1,i-1}$.
Thus $P_{R_1,i}$ are distinct for all $i$.

The next step is to construct a sequence of monoidal transforms of $R$ along $V$ such that the induced chain of prime ideals are not only distinct, but also regular. First notice that if $R_1\rightarrow R_2$ is a single monoidal transform along $V$, then $(V/P_{V,i})_{P_{V,i}}$ is algebraic over $(R_2/P_{R_2,i})_{P_{R_2,i}}$ for all $i$, since $(R_2/P_{R_2,i})_{P_{R_2,i}}$ is an intermediate field between $(R_1/P_{R_1,i})_{P_{R_1,i}}$ and
$(V/P_{V,i})_{P_{V,i}}$. Thus the condition (\ref{thiscondition}) is preserved by further monoidal transforms of
$R_1$ along $V$.

Note that the strict transform of $P_{R_1,i}$ in $R_2$ is $ P_{V,i}\cap R_2=P_{R_2,i}$, if $R_1\rightarrow R_2$ is
centered at a regular prime ideal $\overline{a}$ which properly contains $P_{R_1,i}$.

When $t>1$, we apply Theorem 2.9  to construct a sequence of monoidal transforms $R_1 \rightarrow R_2$ along $V$, centered at ideals which properly contain the strict transform of  $P_{R_1,t-1}$, such that the strict transform of  $P_{R_1,t-1}$in $R_2$ is a regular prime. This regular prime is necessarily $P_{R_2,t-1}$.

We apply this argument for $P_{R_2,t-2}$ and recursively for all $t-i$ up to $i=t-1$ to construct  sequences of monoidal  transforms
$$
R \rightarrow R_1 \rightarrow \dots \rightarrow R_t$$
along $V$ such that $(V/P_{V,i})_{P_{V,i}}$ is algebraic over $ (R_t/P_{R_t,i})_{P_{R_t,i}}$ for all $i$, and  $P_{R_t,i}$ are regular and distinct prime ideals in $R_t$ for all $i$.

Now we will assume that $R$ has properties (\ref{algebraic}) and (\ref{regdistinct}) since otherwise, we can replace $R$ with $R_t$.
We will proceed to prove the Theorem by induction on $t= \text{rank } \nu$.

The proof when $t=1$ follows from Theorem 6.5 \cite{CG} in the case when $K=K^*$ and $R^*=S^*$.

Assume that the Theorem is true for valuation rings  of rank less than $t$.
$V_{P_{V,1}}$ is a valuation ring of $K$ of rank $t-1$ (since its value group is $\Gamma_t/\Gamma_1$), which dominates $R_{P_{R,1}}$.
By the induction statement, there exists a sequence of monoidal transforms
$R_{P_{R,1}} \rightarrow T_1$
along $V_{P_{V,1}}$ such that the conclusions of the Theorem hold for $T_1$ and $V_{P_{V,1}}$.
Let $P_{V_{P_{V,1}},j}=P_{V,j+1}V_{P_{V,1}}$ for $j=0, \dots, t-1$. We have $P_{V_{P_{V,1}},0}=m_{V_{P_{V,1}}}$.
The induction statement of the Theorem gives us that
$Q(\widehat{(T_1)_{P_{T_1,j}}})\subset \widehat{(T_1)_{P_{T_1,i}}}$ are regular primes for all $j=0, \dots , t-1$.

By Theorem 2.9 \cite{C}, there exists a sequence of monoidal transforms $R \rightarrow R_1$ along $V$
such that $(R_1)_{P_{R_1,1}}\cong T_1$ and $R_1/P_{R_1,i}$ are regular local rings for all $i$.
Thus
$(T_1)_{P_{T_1,j}}\cong (R_1)_{P_{R_1,j+1}}$ and $Q(\widehat{(R_1)_{P_{R_1,j+1}}})$ are regular primes for the valuation ring $(V_{P_{V,1}})_{P_{V,j+1}}V_{P_{V,1}}\cong V_{P_{V,j+1}}$ for $j=0, \dots, t-1$.

Let $A_1=R_1/P_{R_1,1}$. This ring is dominated by the rank $1$ valuation ring $V/P_{V,1}$.
By Theorem 6.5 \cite{CG}  with $K=K^*$ and $R^*=S^*$ there exists a sequence of monoidal transforms along $V/P_{V,1}$
\begin{equation}\label{seq}
A_1\rightarrow A_2 \rightarrow \dots \rightarrow A_r
\end{equation}
such that $A_r$ is a regular local ring and the conclusions of the Theorem hold for $A_r$ and $V/P_{V,1}$.

Now we will show that there exists a sequence of monoidal transforms along $V$
$$R_1\rightarrow R_2 \rightarrow \dots \rightarrow R_r$$
such that $A_r \cong R_r/P_{R_r,1}$ and $(R_r)_{P_{R_r,1}}\cong T_1$.

For $i=1, \dots, r-1$, $A_{i+1}$ is a local ring of the blow-up of a non-zero regular prime ideal $a_i \subset A_i$.

Let $i=1$. Let $\overline{a_1}$ be the preimage of $a_1$ in $R_1$.
The ideal $\overline{a_1}$ is a regular prime in $R_1$ (since  $R_1/\overline{a_1}\cong A_1/a_1$)
which properly contains $P_{R_1,1}$.

 Choose $\overline{a_1}$ to be the center of the blow-up of $R_1$ along $V$. Let $R_2$ be the local ring of this blow-up which is dominated by $V$. We have that $R_2$ is a regular local ring, $P_{R_2,i}$ is the strict transform of $P_{R_1,i}$ in $R_2$ for $i=1,\dots , t$, and  thus $P_{R_2,i}$ are regular primes for $i=0,\dots , t$.
Moreover, by the universal property of blowing up,
 $A_2\cong R_2/P_{R_2,1}$. Further, since  $T_1\cong (R_1)_{P_{R_1,1}}$ and $P_{R_1,1}$ is a proper subset of $\overline {a_1}$, $(R_2)_{P_{R_2,1}}\cong (R_1)_{P_{R_1,1}}\cong T_1$.

Now, we repeat this construction along (\ref{seq}), and obtain a sequence of monoidal transforms
$$
R_1\rightarrow R_2 \rightarrow \dots \rightarrow R_r$$ along $V$ such that
$$
A_i\cong R_i/P_{R_i,1} \text{ and } (R_i)_{P_{R_i,1}}\cong T_1 \text{ for } i=1, \dots , r.
$$
Thus the conclusions of the theorem hold for $R_r$, since $A_r\cong R_r/P_{R_r,1}$ and $\widehat{R_r}/Q(\widehat{R_r})\cong \widehat{A_r}/Q(\widehat{A_r})$.
\end{proof}

\end{document}